\documentclass[12pt]{article}
\newcommand{\be}{\begin{equation}}
\newcommand{\ee}{\end{equation}}
\begin{document}
\title{\textbf {Hurwitz - Bernoulli Numbers, Formal Groups and the L - Functions of Elliptic Curves}}
\author{ H. Gopalakrishna Gadiyar and R. Padma\\
School of Advanced Sciences,\\ V.I.T. University, Vellore 632014 INDIA\\E-mail: gadiyar@vit.ac.in, rpadma@vit.ac.in}
\date{~~}
\maketitle

\begin{abstract}
Classically, Euler developed the theory of the Riemann zeta - function using as his starting point the exponential and partial fraction forms of $\cot z $. In this paper we wish to develop the theory of $L$-functions of elliptic curves starting from the theory of elliptic functions in an analogous manner.

\end{abstract}

The importance of the $L$-function of elliptic curves is well known because of its connection to the conjectures of Shimura,Taniyama and Weil. The centrality of the Bernoulli numbers to many branches of mathematics has been explained in \cite{Mazur}. In this brief note we wish to state the fact the arithmetic of elliptic curves can be developed in analogy with the theory of the Riemann zeta - function. The result we wish to state is that if we take the formal exponential 
\be
\frac{T}{f_E(T)} = \sum_{n=0}^\infty \frac{BH_n}{n!} T^n
\ee
which arises from the elliptic function corresponding to an elliptic curve with the rational $g_2$ and $g_3$, the $L$-function is related to the formal logarithm by the formula 
\be
f_L(T)=\sum_{n=1}^\infty \frac{a_n}{n}T^n
\ee
where, $f_E(f_L(T))=T$, $f_L(f_E(T)=T$ and $L_E(s)={\displaystyle \sum_{n=1}^\infty \frac{a_n}{n^s}}$.

It is well known that the classical theory of elliptic functions was developed in strict analogy with the theory of trigonometric functions. This is evident from the well known formulae we reproduce below \cite{Whittaker:Watson}.
\begin{equation}
\sin (z) = z \prod_{n=-\infty}^{\infty ~\prime}\left ( 1-\frac{z}{m\pi}\right ) e^{\left(\frac{z}{m\pi}\right)} \, ,  
\end{equation}
\begin{equation}
\sigma (z) = z \prod_{m,n}^{\prime} \left(1-\frac{z}{\Omega_{m,n}}\right) exp\left(\frac{z}{\Omega_{m,n}}+\frac{z^2}{2\Omega_{m,n}^2}\right)
\end{equation}
where $\Omega_{m,n}=2m\omega_1+2n\omega_2$, $\omega_1$ and $\omega_2$ are the fundamental periods of the elliptic function. 
Corresponding to the exponential form $\sin (z) = \frac{(e^{iz}-e^{-iz}) }{2i}$, we have the exponential form for $\sigma (z)$
\begin{equation}
\sigma (z)=\frac{2\omega_1}{\pi} exp\left(\frac{\eta_1 z^2}{2\omega_1}\right) \sin \left(\frac{\pi z}{2\omega_1}\right)\prod_{n=1}^\infty \left\{ \frac{1-2q^{2n}\cos \frac{\pi z}{\omega_1}+q^{4n}}{(1-q^{2n})^2}\right\}
\end{equation}
where $q=exp(\frac{i\pi \omega_2}{\omega_1})$ and $\eta_1 = \zeta (\frac{\omega_1}{2})$.
Similarly,
\begin{equation}
\cot (z)=\frac{1}{z}+\sum_{m=-\infty}^{\infty ~\prime}\frac{1}{z-m\pi}+\frac{1}{m\pi}
\end{equation}
and
\be
\zeta (z)= \frac{1}{z}+\sum_{m,n}^{\prime} \left\{ \frac{1}{z-\Omega_{m,n}}+\frac{1}{\Omega_{m,n}}+\frac{z}{\Omega_{m,n}^2}\right\}
\ee
and finally,
\be
\mathrm{cosec} ^2(z)=\frac{1}{z^2}+\sum_{m=-\infty}^{\infty ~\prime} \frac{1}{(z-m\pi)^2}
\ee
\be
\wp (z)=\frac{1}{z^2}+\sum_{m,n}^{~\prime} \left\{ \frac{1}{(z-\Omega_{m,n})^2}-\frac{1}{\Omega_{m,n}^2}\right\}
\ee

In the theory of formal groups \cite{Hazewinkel}, the formal exponential and the formal logarithm have very important properties. In the paper of A. Baker \cite{Baker} these properties have been elucidated in a very general setting. There is also a theorem of Honda \cite{Honda} and \cite{Cartier}, which connects the $L$-function of an elliptic curve to the formal logarithm ${\displaystyle \sum_{n=1}^\infty \frac{a_n}{n}T^n}$. Further, given a formal group law $F(X,Y)$, it is known that 
\be
F(X,Y)=f_E(f_L(X)+f_L(Y)) 
\ee
We assert that the natural way of understanding the arithmetic of elliptic curves is to take the formal exponential to be the one which arises from the Weierstrass elliptic functions. It is well known \cite{Clarke}, \cite{Hazewinkel} and \cite{Silverman} that the classical zeta function arises from 
\be
\frac{T}{f_E(T)} = \sum_{n=0}^\infty \frac{B_n}{n!} T^n
\ee
where $f_E(T) = e^T-1$ and the formal logarithm is given by $\log (1+T)$. This formal group law arises from the trigonometric function $\cot z$. In analogy, it would seem natural that the object in the case of elliptic functions would be got by substituting Bernoulli numbers by Hurwitz-Bernoulli numbers \cite{Hurwitz}, \cite{Katz} and \cite{Lemmermeyer}. That is,
\be
\frac{T}{f_E(T)} = \sum_{n=0}^\infty \frac{BH_n}{n!} T^n
\ee
the formal logarithm would be ${\displaystyle \sum_{n=1}^\infty \frac{a_n}{n}T^n}$. Hurwitz - Bernoulli numbers have been defined in a general setting by Katz\cite{Katz} as given below. For the elliptic curve
\be
y^2= 4x^3-g_2x-g3 \, ,
\ee
we have
\be
\wp (z) = \frac{1}{z^2}+2 \sum_{n\ge 1} G_{2n+2} \frac{z^{2n}}{2n!} \, ,
\ee
and
\be
BH_k = 2kG_k {\rm ~for~} k \ge 4; =0 {\rm ~if~} k {\rm ~is~odd.}
\ee
This would give direct connection between the Weierstrass elliptic functions and the $a_n$'s which trap the arithmetic properties of the elliptic curves via the Hurwitz - Bernoulli numbers. 

It is possible to envisage a general frame work where other conjectures related to the arithmetic of elliptic curves will find a natural setting in the theory of formal groups. See \cite{Hazewinkel} for the discussion of the Birch - Swinnerton -Dyer conjecture. 
\bigskip

\noindent {\bf Acknowledgements} We would like to thank Dr. Jayandra Bandyopadhyay for making available copies of earlier papers on Hurwitz - Bernoulli numbers referred to in \cite{Lemmermeyer}. We would also like to thank Professors Arul Lakshminarayanan (IIT, Madras), A. Sankaranarayanan (TIFR) ,K. Srinivas (IMSc, Chennai) and Dr. Vinosh Babu James for making available the papers and books used in this research.

\end{document}